\documentclass[12pt]{article}
\usepackage{lineno,hyperref}
\modulolinenumbers[5]
\usepackage{color}

\usepackage{graphics,amsmath,amssymb}
\usepackage{amsthm}
\usepackage{amsfonts}
\usepackage{latexsym}
\usepackage{epsf}

\newcommand{\sstirling}[2]{\genfrac\{\}{0pt}{}{#1}{#2}}

\begin{document}

\begin{center}
\vskip 1cm{\LARGE\bf Spivey's Bell Number Formula Revisited}
\vskip 1cm
\large 
Mahid M. Mangontarum\\
Department of Mathematics\\
Mindanao State University---Main Campus\\
Marawi City 9700\\
Philippines \\
\href{mailto:mmangontarum@yahoo.com}{\tt mmangontarum@yahoo.com} \\
\href{mailto:mangontarum.mahid@msumain.edu.ph}{\tt mangontarum.mahid@msumain.edu.ph} \\
\end{center}

\vskip .2 in

\begin{abstract}
This paper introduces an alternative form of the derivation of Spivey's
Bell number formula, which involves the $q$-Boson operators $a$ and
$a^{\dagger}$. Furthermore, a similar formula for the case of the
$(q,r)$-Dowling polynomials is obtained, and is shown to produce a
generalization of the latter.
 \end{abstract}

\section{Introduction}

Consider the Stirling numbers of the second kind, denoted by $\sstirling{m}{j}$, which appear as coefficients in the expansion of
\begin{equation}
t^n=\sum_{k=0}^n\sstirling{n}{k}(t)_k,
\end{equation}
where $(t)_k=t(t-1)(t-2)\cdots(t-k+1)$.The Bell numbers, denoted by $B_n$, are defined by
\begin{equation}
B_n=\sum_{j=0}^{n}\sstirling{n}{j}\label{Belldef}
\end{equation}
and are known to satisfy the recurrence relation
\begin{equation}
B_{n+1}=\sum_{k=0}^n\binom{n}{k}B_k.\label{Bellrec}
\end{equation}
In 2008, Spivey \cite{Spivey} obtained a remarkable formula which unifies the defining relation in \eqref{Belldef} and the identity \eqref{Bellrec}. The said formula is given by
\begin{equation}
B_{n+m}=\sum_{k=0}^n\sum_{j=0}^mj^{n-k}\sstirling{m}{j}\binom{n}{k}B_k\label{spivey}
\end{equation}
and is popularly known as ``Spivey's Bell number formula''. Equation \eqref{spivey} was proved in \cite{Spivey} using a combinatorial approach involving partition of sets. Different proofs and extensions of \eqref{spivey} were later on studied by several authors. For instance, a proof which made use of generating functions was done by Gould and Quaintance \cite{Gould} which was then generalized by Xu \cite{Xu} using Hsu and Shuie's \cite{Hsu} generalized Stirling numbers. Belbachir and Mihoubi \cite{Belbachir} presented a proof that involves decomposition of the Bell polynomials into a certain polynomial basis. Mez\H{o} \cite{Mez} obtained a generalization of the Spivey's formula in terms of the $r$-Bell polynomials via combinatorial approach. The notion of dual of \eqref{spivey} was also presented in the same paper. On the other hand, the work of Katriel \cite{Katriel} involved the use of the operator $X$ satisfying
\begin{equation}
DX-qXD=1,
\end{equation}
where $D$ is the $q$-derivative defined by
\begin{equation}
Df(x)=\frac{f(qx)-f(x)}{x(q-1)}.
\end{equation}
For the sake of clarity and brevity, this method will be referred to as ``Katriel's proof''.

Now, aside from being implicitly implied in Katriel's proof, none of the previously-mentioned studies considered establishing $q$-analogues. It is, henceforth, the main purpose of this paper to obtain a generalized $q$-analogue of Spivey's Bell number formula. 

\section{Alternative form of ``Katriel's proof''}

We direct our attention to the $q$-Boson operators $a$ and $a^\dagger$ satisfying the commutation relation
\begin{equation}
[a,a^\dagger]_q=aa^\dagger-qa^\dagger a=1\label{qBoson}
\end{equation}
(see \cite{Arik}). We define the Fock space (or Fock states) by the basis $\{\left.|s\right\rangle; s=0,1,2,\ldots\}$ so that the relations $a\left.|s\right\rangle=\sqrt{[s]_q}\left.|s-1\right\rangle$ and $a^\dagger\left.|s\right\rangle=\sqrt{[s+1]_q}\left.|s+1\right\rangle$ form a representation that satisfies \eqref{qBoson}. The operators $a^\dagger a$ and $(a^\dagger)^ka^k$, when acting on $\left.|s\right\rangle$, yield
\begin{equation}
a^\dagger a\left.|s\right\rangle=[s]_q\left.|s\right\rangle
\end{equation}
and
\begin{equation}
(a^\dagger)^ka^k\left.|s\right\rangle=[s]_{q,k}\left.|s\right\rangle,\label{qBp1}
\end{equation}
respectively, where $[s]_q=\frac{q^s-1}{q-1}$ and $[s]_{q,k}=[s]_q[s-1]_q[s-2]_q\cdots[s-k+1]_q$. Hence, the $q$-Stirling numbers of the second kind $\sstirling{n}{k}_q$ \cite{Car3} can be defined alternatively as
\begin{equation}
(a^\dagger a)^n=\sum_{k=1}^{n}\sstirling{n}{k}_q(a^\dagger)^ka^k.\label{qStirling}
\end{equation}
From \eqref{qBoson}, it is clear that
\begin{equation}
[a,(a^\dagger)^k]_{q^k}=[a,(a^\dagger)^{k-1}]_{q^{k-1}}a^\dagger+q^{k-1}(a^\dagger)^{k-1}[a,a^\dagger]_q,
\end{equation}
and by induction on $k$, we have
\begin{equation}
[a,(a^\dagger)^k]_{q^k}=[k]_q(a^\dagger)^{k-1}.\label{lem1}
\end{equation}
Since $a\left.|0\right\rangle=0$, then by \eqref{lem1},
\begin{eqnarray*}
a(a^\dagger)^\ell\left.|0\right\rangle&=&[a,(a^\dagger)^{\ell}]_{q^{\ell}}\left.|0\right\rangle\\
&=&[\ell]_q(a^\dagger)^{\ell-1}\left.|0\right\rangle.
\end{eqnarray*}
Moreover,
\begin{equation}
a^k(a^\dagger)^\ell\left.|0\right\rangle=\frac{[\ell]_q!}{[\ell-k]_q!}(a^\dagger)^{\ell-k}\left.|0\right\rangle,
\end{equation}
for $k\leq\ell$ and
\begin{equation}
a^k(a^\dagger)^\ell\left.|0\right\rangle=0,
\end{equation}
for $k>\ell$. Finally,
\begin{equation}
a^ke_q(xa^\dagger)\left.|0\right\rangle=x^ke_q(xa^\dagger)\left.|0\right\rangle,\label{lem2}
\end{equation}
where $e_q(xa^\dagger)$ is the $q$-exponential function defined by
\begin{equation}
e_q(t)=\sum_{\ell=0}^\infty\frac{t^\ell}{[\ell]_q!}.
\end{equation}
Applying \eqref{lem2} to \eqref{qStirling} yields
\begin{equation}
(a^{\dagger}a)^ne_q(ta^\dagger)\left.|0\right\rangle=B_{n,q}(ta^\dagger)e_q(ta^\dagger)\left.|0\right\rangle,
\end{equation}
where $B_{n,q}(ta^{\dagger})$ denotes the $q$-Bell polynomials defined by
\begin{equation}
B_{n,q}(t)=\sum_{k=0}^n\sstirling{n}{k}_qt^k.
\end{equation}
Let $x=ta^\dagger$ so that
\begin{equation}
(a^{\dagger}a)^ne_q(x)\left.|0\right\rangle=B_{n,q}(x)e_q(x)\left.|0\right\rangle.\label{qBp2}
\end{equation}
Before proceeding, note that by definition,
\begin{equation}
[a,(a^\dagger)^k]_{q^k}=a(a^\dagger)^k-q^k(a^\dagger)^ka.
\end{equation}
By \eqref{lem1},
\begin{eqnarray*}
a(a^\dagger)^k-q^k(a^\dagger)^ka&=&[k]_q(a^\dagger)^{k-1}\\
a(a^\dagger)^k&=&q^k(a^\dagger)^ka+[k]_q(a^\dagger)^{k-1}.
\end{eqnarray*}
This can be further expressed as
\begin{equation}
(a^\dagger a)(a^\dagger)^k=(a^\dagger)^k\big([k]_q+q^k(a^\dagger a)\big).\label{lem3}
\end{equation}
Now, we have
\begin{eqnarray*}
(a^{\dagger}a)^{n+m}&=&(a^{\dagger}a)^n\sum_{j=0}^m\sstirling{m}{j}_q(a^{\dagger})^ja^j\\
&=&\sum_{j=0}^m\sstirling{m}{j}_q(a^{\dagger})^j\big([j]_q+q^j(a^\dagger a)\big)^na^j\\
&=&\sum_{j=0}^m\sum_{k=0}^n\sstirling{m}{j}_q\binom{n}{k}[j]_q^{n-k}q^{jk}(a^{\dagger})^j(a^\dagger a)^ka^j.
\end{eqnarray*}
Multiplying both sides with $e_q(x)\left.|0\right\rangle$ makes the left-hand side
\begin{equation}
(a^{\dagger}a)^{n+m}e_q(x)\left.|0\right\rangle=B_{n+m,q}(x)e_q(x)\left.|0\right\rangle,
\end{equation}
while the right-hand side becomes
\begin{eqnarray*}
\sum_{j=0}^m\sum_{k=0}^n\sstirling{m}{j}_q\binom{n}{k}[j]_q^{n-k}q^{jk}(a^{\dagger})^j(a^\dagger a)^ke_q(x)\left.|0\right\rangle a^j&=&\sum_{j=0}^m\sum_{k=0}^n\sstirling{m}{j}_q\binom{n}{k}[j]_q^{n-k}q^{jk}\\
& &B_{k,q}(x)e_q(x)\left.|0\right\rangle(a^{\dagger})^ja^j.
\end{eqnarray*}
Dividing both sides by $e_q(x)\left.|0\right\rangle$ and using \eqref{qBp1} gives
\begin{equation}
B_{n+m,q}(x)=\sum_{j=0}^m\sum_{k=0}^n\sstirling{m}{j}_q\binom{n}{k}[j]_q^{n-k}q^{jk}B_{k,q}(x)[x]_{q,j}.\label{result1}
\end{equation}
As $q\rightarrow 1$, we obtain a polynomial version of Spivey's Bell number formula which, in return, reduces to \eqref{spivey} when we set $x=1$.

It is important to emphasize that this is not a new proof, but an
alternative form of Katriel's proof, since the operators $a$,
$a^{\dagger}$ and the operators $X$, $D$ generate isomorphic algebras.

\section{A generalization of Spivey's Bell number formula}

The main result of this paper is the following identity:
\begin{equation}
D_{m,r,q}(n+\ell,x)=\sum_{j=0}^\ell\sum_{k=0}^nm^jW_{m,r,q}(\ell,j)\binom{n}{k}(m[j]_q+r)^{n-k}q^{jk}D_{m,0,q}(k,x)[x]_{q,j}.\label{result2}
\end{equation}
Here, $D_{m,r,q}(n,x)$ is a $(q,r)$-Dowling polynomial
defined previously by the author and Katriel \cite{MangontarumKatriel} as
\begin{equation}
D_{m,r,q}(n,x)=\sum_{k=0}^nW_{m,r,q}(n,k)x^k,
\end{equation}
where $W_{m,r,q}(n,k)$ is the $(q,r)$-Whitney numbers of the second kind. Several properties of $D_{m,r,q}(n,x)$ can be seen in \cite{Mangontarum,MangontarumKatriel}. 

To derive \eqref{result2}, we first multiply both sides of \eqref{lem3} by $m$ and then add $r(a^{\dagger})^k$ to yield  
\begin{equation}
(ma^\dagger a+r)(a^\dagger)^k=(a^\dagger)^k(m[k]_q+r+mq^ka^\dagger a).\label{lem4}
\end{equation}
Also, multiplying both sides of the defining relation in \cite[Equation\ 16]{MangontarumKatriel} by $e_q(ta^{\dagger})\left.|0\right\rangle$ and applying \eqref{lem2} yields
\begin{eqnarray*}
(ma^\dagger a+r)^ne_q(ta^\dagger)\left.|0\right\rangle&=&\sum_{k=0}^nm^kW_{m,r,q}(n,k)(a^\dagger)^k a^ke_q(ta^\dagger)\left.|0\right\rangle\\
&=&\sum_{k=0}^nm^kW_{m,r,q}(n,k)(a^\dagger)^kt^ke_q(ta^\dagger)\left.|0\right\rangle\\
&=&D_{m,r,q}(n,mta^\dagger)e_q(ta^\dagger)\left.|0\right\rangle.
\end{eqnarray*}
Now, by \eqref{lem4},
\begin{eqnarray*}
(ma^\dagger a+r)^{n+\ell}&=&\sum_{j=0}^\ell m^jW_{m,r,q}(\ell,j)(ma^\dagger a+r)^n(a^\dagger)^ja^j\\
&=&\sum_{j=0}^\ell m^jW_{m,r,q}(\ell,j)(a^\dagger)^j(m[j]_q+r+mq^ja^\dagger a)^na^j\\
&=&\sum_{j=0}^\ell\sum_{k=0}^nm^{j+k}W_{m,r,q}(\ell,j)\binom{n}{k}(a^\dagger)^j(m[j]_q+r)^{n-k}q^{kj}(a^\dagger a)^ka^j.
\end{eqnarray*}
Applying this expression to the operator identity $e_q(ta^{\dagger})\left.|0\right\rangle$, combining with the previous equation, using \eqref{qBp1}, \eqref{qBp2} and $W_{m,0,q}(k,i)=m^{k-i}\sstirling{k}{i}_q$ (see \cite[Equation\ 18]{MangontarumKatriel}), and then dividing both sides of the resulting identity by $e_q(ta^{\dagger})\left.|0\right\rangle$ completes the derivation.

\section{Remarks}

Since $W_{1,0,q}(\ell,j)=\sstirling{\ell}{j}_q$, then by setting $x=1$, $m=1$ and $r=0$, we have
\begin{equation}
D_{1,0,q}(n+\ell,1)=\sum_{j=0}^\ell\sum_{k=0}^n\sstirling{\ell}{j}_q\binom{n}{k}[j]_q^{n-k}q^{jk}B_{k,q},
\end{equation}
where $B_{k,q}:=B_{k,q}(1)$. This is a $q$-analogue of \eqref{spivey} which was first obtained by Katriel \cite{Katriel}. On the other hand, setting $x=1$ and then taking the limit of \eqref{result2} as $q\rightarrow 1$ provides a generalization of Spivey's Bell number formula in terms of the $r$-Whitney numbers of the second kind, denoted by $W_{m,r}(\ell,j)$, and the $r$-Dowling numbers, denoted by $D_{m,r}(n)$, (see \cite{Cheon,Mez1}), given by
\begin{equation}
D_{m,r}(n+\ell)=\sum_{j=0}^\ell\sum_{k=0}^nm^jW_{m,r}(\ell,j)\binom{n}{k}(mj+r)^{n-k}D_{m,0}(k).\label{result3}
\end{equation}

In a recent paper, Mansour et al. \cite{Mansour} obtained the following generalization of Spivey's Bell number formula:
\begin{equation}
D_{p,q}(a+b;x)=\sum_{i=0}^{a}\sum_{j=0}^{b}\sum_{\ell=0}^{j}(mq^i)^{j-\ell}x^{i+\ell}\binom{b}{j}\left([r]_p+m[i]_q\right)^{b-j}W_{p,q}(a,i)S_q(j,\ell).\label{pqMansour}
\end{equation}
Here, $D_{p,q}(n;x)$ and $W_{p,q}(n,k)$ denote the $(p,q)$-analogues of the $r$-Dowling polynomials and the $r$-Whitney numbers of the second kind, respectively. The $(p,q)$-analogues are natural generalizations of $q$-analogues. However, since the manner by which the numbers $W_{m,r,q}(n,k)$ were defined in \cite{MangontarumKatriel} differs from the work of Mansour et al. \cite{Mansour}, the main result of this paper is not generalized by \eqref{pqMansour}.

\section{Acknowledgment}
The author is very thankful to the editor-in-chief and to the referee(s) for carefully reading the paper. Their comments and suggestions were very helpful. Special thanks also to Dr.~Jacob Katriel for his insights on $q$-Boson operators. This paper is dedicated to my family and all other victims of the war in Marawi.

\end{document}